\documentclass{amsart}

\usepackage{latexsym,enumerate}
\usepackage{amsmath,amsthm,amsopn,amstext,amscd,amsfonts,amssymb}
\usepackage{verbatim}
\usepackage{graphicx}
\usepackage{pstricks}
\usepackage{wasysym}
\usepackage[all]{xy}
\usepackage{epsfig} 
\usepackage{graphicx,psfrag} 
\usepackage{subfigure}
\usepackage{pstricks,pst-node}
\usepackage{euscript}
\usepackage{amsmath}
\numberwithin{equation}{section}
\usepackage{amsfonts}
\usepackage{amssymb}
\usepackage{graphicx}
\usepackage{amscd}
\usepackage{commath}
\usepackage{tikz}
\usepackage{mathtools}
\usepackage{physics}
\usepackage{hyperref}

\usepackage{multicol}
\usepackage{color}
\usepackage{colortbl}

\newcommand{\NN}{\mathbb{N}}

\newcommand{\RR}{\mathbb{R}}
\newcommand{\ZZ}{\mathbb{Z}}

\newtheorem{theorem}{Theorem}[section]
\newtheorem{lemma}[theorem]{Lemma}
\newtheorem{proposition}[theorem]{Proposition}

\newtheorem{definition}[theorem]{Definition}
\newtheorem{example}[theorem]{Example}
\newtheorem{obs}[theorem]{Remark}

\newcommand{\spb}[1]{\smallskip}
\newcommand{\mpb}[1]{\medskip}
\newcommand{\bpb}[1]{\bigskip}

\begin{document}
\DeclareGraphicsExtensions{.jpg,.pdf,.mps,.png}

\title{Gluing and decomposition of doubling metric graphs and metric spaces}

\author[Alvaro Mart\'{\i}nez-P\'erez]{Alvaro Mart\'{\i}nez-P\'erez$^{(1)}$}
\address{Departamento de Álgebra Geometría y Topología, Universidad Complutense de Madrid, Plaza de las Ciencias 3, 28040, Madrid
Spain. ORCID: 0000-0002-1344-6189.}
\email{alvaro36@ucm.es}
\thanks{$^{(1)}$ Supported in part by a grant
from the Agencia Estatal de Investigaci\'on (PID2024-156663NB-I00), Spain.
}

\author[Amaia Zudaire]{Amaia Zudaire$^{(1)}$}
\address{Facultad de Ciencias Matemáticas, Universidad Complutense de Madrid, Plaza de las Ciencias 3, 28040, Madrid
Spain.}
\email{amaiazud@ucm.es}
\thanks{}

\date{\today}

\begin{abstract}
Herein we study the conditions under which the union (or gluing) of a family of (metric) doubling graphs yields a doubling graph. Conversely, we consider a doubling graph obtained by gluing a family of graphs and find conditions under which those graphs must also be doubling. Finally, we extend these results to the more general setting of metric spaces with an appropriate definition of gluing. 
\end{abstract}

\maketitle{}



\section{Introduction}

A metric space $(X,d)$ is said to be \emph{(metric) doubling} if there exists
a constant $C \geq 1$ such that every ball of radius $2r$ can be
covered by at most $C$ balls of radius $r$, for every $r>0$. If a doubling measure exists on a metric space $X$, then $X$ itself is a metric doubling space. Conversely, every complete metric space admits a doubling measure. 
The doubling condition
expresses a quantitative form of finite dimensionality: at every
location and scale, the number of balls needed to cover a larger
ball is uniformly controlled. This property serves as a geometric analogue to spaces with bounded finite dimensionality, capturing the idea that the volume of local neighborhoods does not grow exponentially with respect to the radius. The infimum of $\log_2 C$ over all admissible constants $C$ is
called the \emph{doubling dimension} of $X$. Also, doubling spaces are exactly the spaces with a finite Assouad dimension.  For a wide introduction to this property, we refer to \cite{CW} and \cite{Hei}. It is straightforward to verify that the definition of a doubling space remains equivalent whether it is formulated using open or closed balls. For convenience, throughout this paper, any reference to a ball will implicitly mean a closed ball.

An important class of examples is provided by graphs equipped
with their shortest-path metric. Given a connected graph
$G=(V,E)$, let $\delta_G(x,y)$
be the minimum number of edges in a path joining $x,y \in V$. See \cite{Diestel} for standard terminology and background.
We say that $G$ is a \emph{doubling graph} if the metric space
$(V,\delta_G)$ is doubling. Every finite connected graph is
doubling, although its doubling constant need not be bounded
independently of its size. A standard infinite example is the
integer lattice $\mathbb{Z}^d$ with its usual nearest-neighbor
graph structure, whose shortest-path metric is the $\ell^1$
metric.

Doubling spaces play an important role in metric geometry,
harmonic analysis, and theoretical computer science. See, for example, \cite{Chan,Gill}. Their
uniform control of covering numbers provides a useful framework
for studying geometric and analytic questions across different
scales. In \cite{MPR26}, the authors prove that a uniform graph $G$ satisfies the (Cheeger) isoperimetric inequality, then $G$ is not (metric) doubling and see that the converse is not true. Also, they prove that if $G$ is a doubling graph with doubling constant $C$, then it is $p$-parabolic for every $p \geq \log_2(C)$ and see that the converse is not true. See also \cite{CSC,ppara}.


The aim of this paper is to study conditions under which gluing
doubling spaces preserves the doubling property. Our motivation
comes from previous work on the behavior of geometric and
analytic properties under unions and decompositions of graphs
and metric spaces. In \cite{RT04, PRT04}, the authors study
Gromov hyperbolicity through decompositions of metric spaces,
giving conditions under which the hyperbolicity of the whole
space can be deduced from that of its components.
In \cite{MPR24}, the authors study
$p$-parabolicity on graphs and obtain characterizations in terms
of graph decompositions. They also study its behavior under
vertex identifications.

Following this line of research, we consider the doubling
property, beginning with unions of graphs. 
%
In the graph setting, we use the term union rather than
gluing. Here, the union is understood as follows. Given a pair of graphs $G_1=(V_1,E_1)$ and $G_2=(V_2,E_2)$,  with possibly $V_1\cap V_2\neq \emptyset$, $G:=G_1\cup G_2$ if $G=(V_1\cup V_2,E_1\cup E_2)$. This is equivalent as considering the gluing of two graphs obtained
by identifying the corresponding vertices through a bijection between subsets of $V_1$ and $V_2$. 

The initial focus on graphs is motivated by the flexibility of
this construction. As long as the resulting graph is connected,
its vertex set carries a well-defined shortest-path metric. If the graph is not connected, the metric can be adapted defining an infinite distance between every pair of disconnected vertices. 
Thus, the construction itself does not require the intrinsic
metrics of the original graphs to agree along the identified
parts. This allows us to study the doubling property under
unions that may alter distances within each
constituent graph.

This flexibility also creates the main geometric difficulty.
If a graph $G$ is obtained by the union of $G_1$ and $G_2$, distances
in $G$ between vertices coming from $G_i$ may be much smaller
than their original distances in $G_i$. Indeed, a path in the
resulting graph may leave one constituent graph, travel through the other, and return, creating shortcuts that were not previously
available. Consequently, coverings defined using the intrinsic
metrics of the pieces do not automatically provide suitable
coverings for the glued metric. Controlling the interaction
between the union and these intrinsic metrics is therefore
essential.

When taking the union of infinitely many graphs, the resulting
graph may fail to be doubling, even if the original graphs
have a common doubling constant. We illustrate this with
the following example.

\begin{example}
Let $\{G_i\}_{i \in I}$ be a countably infinite family of graphs,
each consisting of three vertices and two edges joining one
vertex to the other two. These graphs have a common doubling
constant. By identifying their vertices appropriately, one can
construct the infinite rooted binary tree, also known 
as the Cantor tree; see Figure~\ref{can}.

This tree is clearly not doubling. Thus, gluing an
infinite family of uniformly doubling graphs need not produce
a doubling graph.

\begin{figure}[h]
    \centering
    \input{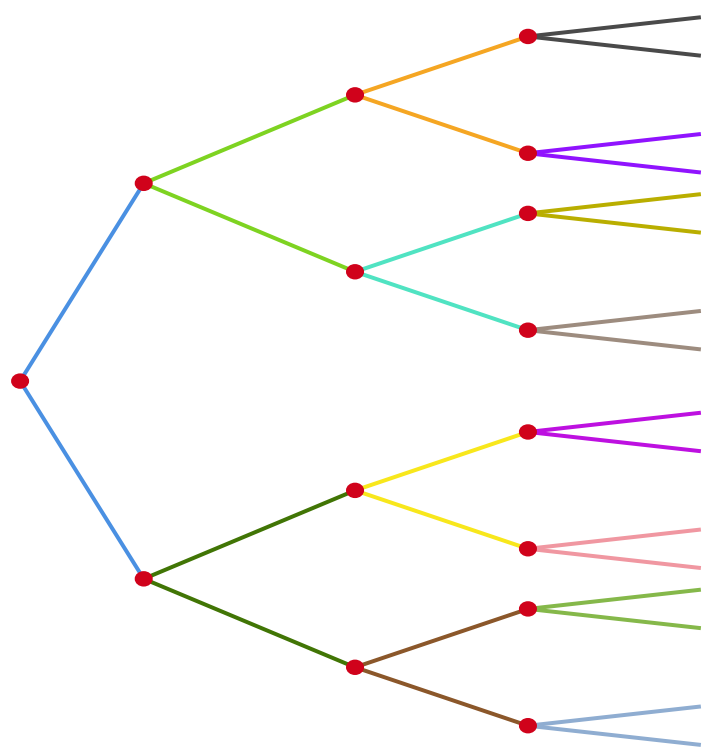}
    \caption{Cantor tree.}
    \label{can}
\end{figure}
\end{example}

For the union of a finite family of doubling graphs, Proposition \ref{p0} and theorems 
\ref{p1} and \ref{p2} provide sufficient conditions ensuring
that the resulting graph remains doubling. These results form
the core of the first two main sections of the paper. We also provide Example \ref{exp} to illustrate that these conditions are not necessary.

In section \ref{s:decomposition}, we consider a doubling graph obtained by gluing a family of graphs and see whether or not these graphs must also be doubling. We provide an example to show that the union of non-doubling graphs can produce a doubling graph and, in Theorem \ref{t:decomp}, we find a sufficient condition so that if the union of an arbitrary family of graphs is doubling, each of the graphs is doubling.

For general metric spaces, gluing requires more care.
After identifying subsets of two metric spaces, the distance
defined by taking infima of lengths of chains may only be
a pseudometric: two distinct points in the glued space may
have distance zero. We therefore need suitable assumptions
to ensure that this construction gives a metric.

In the final section, we address this issue within the framework
of $L$-bi-Lipschitz gluing. The subsets to be identified are
related by a bijection $\varphi \colon A_1 \to A_2$, where
$A_i \subseteq X_i$, satisfying
$$
    L^{-1}d_1(x,y)
    \leq d_2\bigl(\varphi(x),\varphi(y)\bigr)
    \leq Ld_1(x,y),
    \qquad x,y \in A_1,
$$
for some $L \geq 1$. This quantitative compatibility between
the metrics on the identified subsets provides the framework
for Theorem \ref{th:metric}, our final result on the preservation of doubling property.

\section{Union of doubling graphs with finite intersection}
In this section, we establish that a finite union of doubling graphs is also doubling, provided that their pairwise intersections are finite.

\begin{lemma}\label{l:doubling}
$G$ is a doubling graph if and only if there exists a constant $C>0$ such that for any $l\geq 1$, every closed ball of radius $k$ can be covered by at most $C^{l}$ balls of radius $k/2^l$.
\end{lemma}
\begin{proof}
The right-to-left implication is trivial, as the case $l=1$ is precisely the definition of a doubling graph.

For the converse implication, we proceed by induction on $l$. The base case $l=1$ holds by definition. Assume that the statement is true for $l-1$. By the induction hypothesis, any ball of radius $k$ can be covered by at most $C^{l-1}$ balls of radius $k/2^{l-1}$. Furthermore, by definition, each ball of radius $k/2^{l-1}$ can be covered by at most $C$ balls of radius $k/2^l$. Therefore, any ball of radius $k$ can be covered by at most $C^{l-1}C = C^l$ balls of radius $k/2^l$.

\end{proof}

\begin{proposition}\label{p0}
Let $\{G_i\}_{i=1}^N$ be a finite family of doubling graphs such that, for all $i \neq j$, the intersection $G_i \cap G_j$ is finite. Then the graph
$$
G = \bigcup_{i =1}^N G_i
$$
is doubling.
\end{proposition}

\begin{proof}
The proof is by induction on $N$

Base case, $N=2$.

Let $G_1$ and $G_2$ be two doubling graphs such that $\sharp(G_1\cap G_2)=M$. And let $C_1$ and $C_2$ be their doubling constants. 

We take $B$ a ball of $G=G_1\cup G_2$ with radius $2k$.

\begin{itemize}

\item If $B\cap G_1 \neq\emptyset$ and $B\cap G_2=\emptyset$ then $B$ can be covered by $C_1$ balls of radius $k$.

\item If $B\cap G_1=\emptyset$ and $B\cap G_2 \neq \emptyset$ then $B$ can be covered by $C_2$ balls of radius $k$.

\item If $B\cap G_1\neq \emptyset$ and $B\cap G_2 \neq \emptyset$. We take the set $G_1 \cap G_2\cap B=\{x_j\}_{j\in J}$ and for each of theses $x_j$ we define:
$$
B_{j_1}=\{y\in G_1:\delta_{G_1} (x_j,y)\leq 4k\},
$$
$$
B_{j_2}=\{y\in G_2:\delta_{G_2} (x_j,y)\leq 4k\}.
$$
Then $B_{j_1}$ is a ball centered at $x_j$ with radius $4k$ in the graph $G_1$, and, by Lemma \ref{l:doubling},  it can be covered by $C_1^2$ balls of radius $k$.
 Similarly, $B_{j_2}$ can be covered by $C_2^2 $ balls of radius $k$.

Furthermore, 
$$
B\cap G_1 \subset \cup_{j\in J}B_{j_1},
$$
$$
B\cap G_2 \subset \cup_{j\in J}B_{j_2},
$$
so,
$$
B\subset \cup_{j\in J}(B_{j_1}\cup B_{j_2}).
$$

Hence, $B$ can be covered by at most $M(C_1^2+C_2^2)$ balls of radius $k$.

\end{itemize}

Then, $G$ is a doubling graph.

Assume that the stament holds for some arbitrary positive integer $N-1$, i.e. $G' = \bigcup_{i =1}^{N-1} G_i$ is a doubling graph. Let us call $C$ its doubling  constant.

We must show that $G = \bigcup_{i =1}^N G_i$ is a doubling graph. 

Following the same reasoning as the base case and since $G'$ is a doubling graph then $G = \bigcup_{i =1}^N G_i= G'\cup X_N$ is a doubling graph with doubling constant $C^2+C_N^2$.

\end{proof}

However, even if we consider a finite family of doubling graphs, if their intersection is infinite, the union graph need not be doubling.

\begin{example}
Let $G_1=(V_1,E_1)$ be the graph such that $V_1=\{(x,y)\in \ZZ \times \ZZ: y\leq 0\}$ and  $E_1=\{\big((x,y),(z,t)\big)\in V_1\times V_1: x=z \text{  and  } |y-t|=1 \text{  or  } y=t \text{  and  } |x-z|=1\}$.

Let $G_2=(V_2,E_2)$ be the graph such that $V_2=\cup_{i\in \NN\cup \{0\}}(\{(x_i,y)\in \ZZ \times \ZZ: y\geq 0\}\cup \{(-x_i,y)\in \ZZ \times \ZZ: y\geq 0\}$ where $x_0=0$ and $x_i=x_{i-1}+i$. And $E_2=\{((x,y),(x,z))\in V_2 \times V_2:|y-z|=1\} \cup\{\big((x_i,y),(x_{i\pm 1},z)\big)\in V_2\times V_2\} \cup\{\big((-x_i,y),(-x_{i\pm 1},z)\big)\in V_2\times V_2\}$. We construct the union of both graphs, $G$, by identifying the vertices that lie on the line $\RR\times \{0\}$, see figure \ref{u2d}.

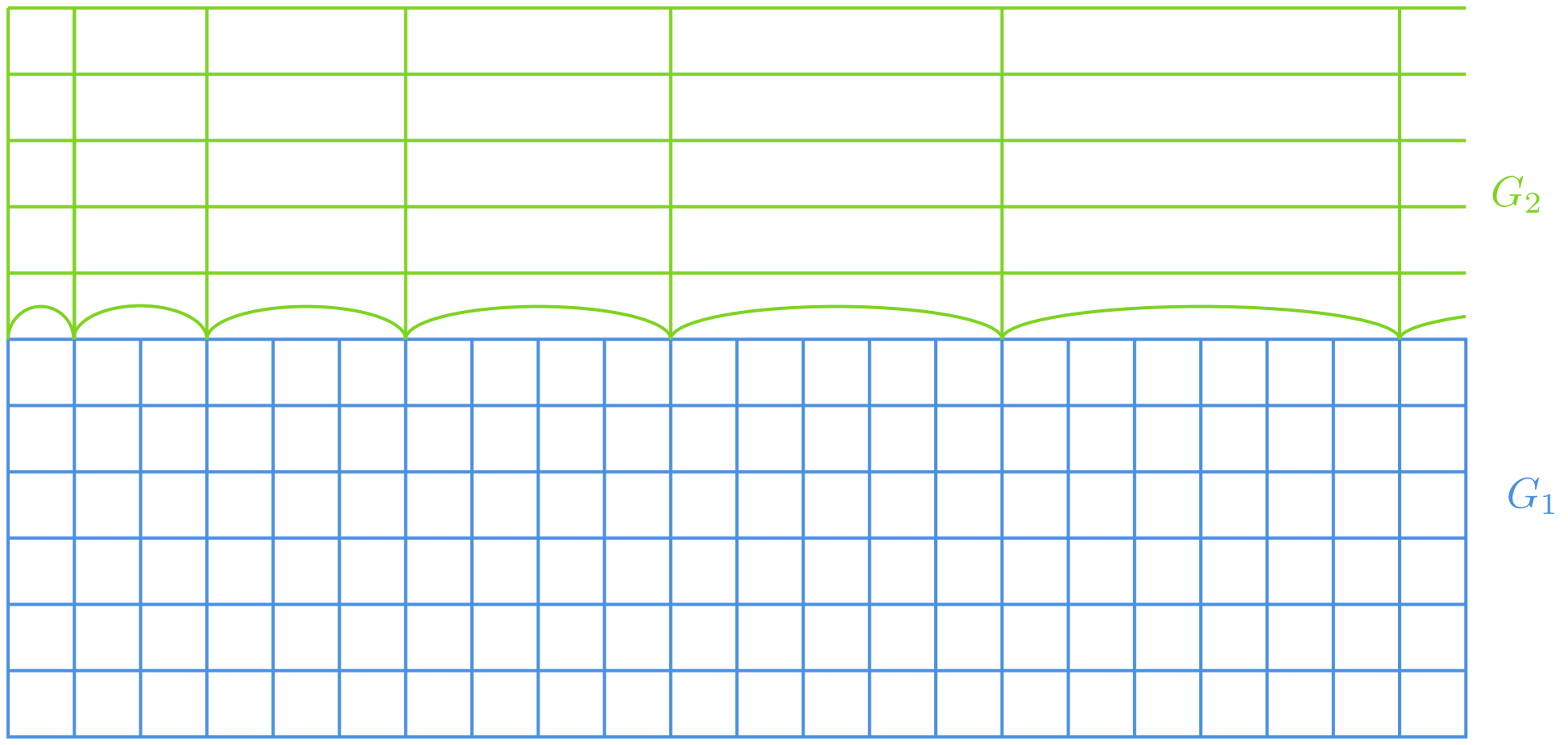
\begin{figure}[h]
\centering
\tikzset{every picture/.style={line width=0.75pt}} 

\begin{tikzpicture}[x=0.75pt,y=0.75pt,yscale=-1,xscale=1]

\draw  [draw opacity=0] (78.28,240.34) -- (78.28,116.46) -- (532.52,116.46) -- (532.52,240.34) -- cycle ; \draw  [color={rgb, 255:red, 74; green, 144; blue, 226 }  ,draw opacity=1 ] (78.28,219.69) -- (532.52,219.69)(78.28,199.05) -- (532.52,199.05)(78.28,178.4) -- (532.52,178.4)(78.28,157.75) -- (532.52,157.75)(78.28,137.1) -- (532.52,137.1) ; \draw  [color={rgb, 255:red, 74; green, 144; blue, 226 }  ,draw opacity=1 ] (98.93,240.34) -- (98.93,116.46)(119.57,240.34) -- (119.57,116.46)(140.22,240.34) -- (140.22,116.46)(160.87,240.34) -- (160.87,116.46)(181.52,240.34) -- (181.52,116.46)(202.16,240.34) -- (202.16,116.46)(222.81,240.34) -- (222.81,116.46)(243.46,240.34) -- (243.46,116.46)(264.1,240.34) -- (264.1,116.46)(284.75,240.34) -- (284.75,116.46)(305.4,240.34) -- (305.4,116.46)(326.05,240.34) -- (326.05,116.46)(346.69,240.34) -- (346.69,116.46)(367.34,240.34) -- (367.34,116.46)(387.99,240.34) -- (387.99,116.46)(408.63,240.34) -- (408.63,116.46)(429.28,240.34) -- (429.28,116.46)(449.93,240.34) -- (449.93,116.46)(470.58,240.34) -- (470.58,116.46)(491.22,240.34) -- (491.22,116.46)(511.87,240.34) -- (511.87,116.46) ; \draw  [color={rgb, 255:red, 74; green, 144; blue, 226 }  ,draw opacity=1 ] (78.28,240.34) -- (78.28,116.46) -- (532.52,116.46) -- (532.52,240.34) -- cycle ;
\draw [color={rgb, 255:red, 126; green, 211; blue, 33 }  ,draw opacity=1 ]   (98.93,116.46) -- (98.93,13.22) ;
\draw [color={rgb, 255:red, 126; green, 211; blue, 33 }  ,draw opacity=1 ]   (98.93,116.46) -- (98.93,13.22) ;
\draw [color={rgb, 255:red, 126; green, 211; blue, 33 }  ,draw opacity=1 ]   (78.28,116.46) -- (78.28,13.22) ;
\draw [color={rgb, 255:red, 126; green, 211; blue, 33 }  ,draw opacity=1 ]   (140.22,116.46) -- (140.22,13.22) ;
\draw [color={rgb, 255:red, 126; green, 211; blue, 33 }  ,draw opacity=1 ]   (202.16,116.46) -- (202.16,13.22) ;
\draw [color={rgb, 255:red, 126; green, 211; blue, 33 }  ,draw opacity=1 ]   (284.75,116.46) -- (284.75,13.22) ;
\draw [color={rgb, 255:red, 126; green, 211; blue, 33 }  ,draw opacity=1 ]   (387.99,116.46) -- (387.99,13.22) ;
\draw [color={rgb, 255:red, 126; green, 211; blue, 33 }  ,draw opacity=1 ]   (511.87,116.46) -- (511.87,13.22) ;
\draw [color={rgb, 255:red, 126; green, 211; blue, 33 }  ,draw opacity=1 ]   (78.28,95.81) -- (532.52,95.81) ;
\draw [color={rgb, 255:red, 126; green, 211; blue, 33 }  ,draw opacity=1 ]   (78.28,75.16) -- (532.52,75.16) ;
\draw [color={rgb, 255:red, 126; green, 211; blue, 33 }  ,draw opacity=1 ]   (78.28,54.52) -- (532.52,54.52) ;
\draw [color={rgb, 255:red, 126; green, 211; blue, 33 }  ,draw opacity=1 ]   (78.28,33.87) -- (532.52,33.87) ;
\draw [color={rgb, 255:red, 126; green, 211; blue, 33 }  ,draw opacity=1 ]   (78.28,13.22) -- (532.52,13.22) ;
\draw  [draw opacity=0] (98.71,116.46) .. controls (98.71,116.46) and (98.71,116.46) .. (98.71,116.46) .. controls (98.71,116.46) and (98.71,116.46) .. (98.71,116.46) .. controls (98.71,110.82) and (94.14,106.24) .. (88.49,106.24) .. controls (82.85,106.24) and (78.28,110.82) .. (78.28,116.46) -- (88.49,116.46) -- cycle ; \draw  [color={rgb, 255:red, 126; green, 211; blue, 33 }  ,draw opacity=1 ] (98.71,116.46) .. controls (98.71,116.46) and (98.71,116.46) .. (98.71,116.46) .. controls (98.71,116.46) and (98.71,116.46) .. (98.71,116.46) .. controls (98.71,110.82) and (94.14,106.24) .. (88.49,106.24) .. controls (82.85,106.24) and (78.28,110.82) .. (78.28,116.46) ;  
\draw  [draw opacity=0] (140.22,116.43) .. controls (140.22,116.43) and (140.22,116.43) .. (140.22,116.43) .. controls (140.22,110.68) and (130.93,106.01) .. (119.46,106.01) .. controls (108,106.01) and (98.71,110.68) .. (98.71,116.43) -- (119.46,116.43) -- cycle ; \draw  [color={rgb, 255:red, 126; green, 211; blue, 33 }  ,draw opacity=1 ] (140.22,116.43) .. controls (140.22,116.43) and (140.22,116.43) .. (140.22,116.43) .. controls (140.22,110.68) and (130.93,106.01) .. (119.46,106.01) .. controls (108,106.01) and (98.71,110.68) .. (98.71,116.43) ;  
\draw  [draw opacity=0] (202.16,116.46) .. controls (202.16,116.46) and (202.16,116.46) .. (202.16,116.46) .. controls (202.16,110.82) and (188.3,106.24) .. (171.19,106.24) .. controls (154.09,106.24) and (140.22,110.82) .. (140.22,116.46) .. controls (140.22,116.46) and (140.22,116.46) .. (140.22,116.46) -- (171.19,116.46) -- cycle ; \draw  [color={rgb, 255:red, 126; green, 211; blue, 33 }  ,draw opacity=1 ] (202.16,116.46) .. controls (202.16,116.46) and (202.16,116.46) .. (202.16,116.46) .. controls (202.16,110.82) and (188.3,106.24) .. (171.19,106.24) .. controls (154.09,106.24) and (140.22,110.82) .. (140.22,116.46) .. controls (140.22,116.46) and (140.22,116.46) .. (140.22,116.46) ;  
\draw  [draw opacity=0] (284.75,116.46) .. controls (284.75,116.46) and (284.75,116.46) .. (284.75,116.46) .. controls (284.75,116.46) and (284.75,116.46) .. (284.75,116.46) .. controls (284.75,110.82) and (266.26,106.24) .. (243.46,106.24) .. controls (220.65,106.24) and (202.16,110.82) .. (202.16,116.46) -- (243.46,116.46) -- cycle ; \draw  [color={rgb, 255:red, 126; green, 211; blue, 33 }  ,draw opacity=1 ] (284.75,116.46) .. controls (284.75,116.46) and (284.75,116.46) .. (284.75,116.46) .. controls (284.75,116.46) and (284.75,116.46) .. (284.75,116.46) .. controls (284.75,110.82) and (266.26,106.24) .. (243.46,106.24) .. controls (220.65,106.24) and (202.16,110.82) .. (202.16,116.46) ;  
\draw  [draw opacity=0] (387.99,116.46) .. controls (387.99,116.46) and (387.99,116.46) .. (387.99,116.46) .. controls (387.99,110.82) and (364.88,106.24) .. (336.37,106.24) .. controls (307.86,106.24) and (284.75,110.82) .. (284.75,116.46) -- (336.37,116.46) -- cycle ; \draw  [color={rgb, 255:red, 126; green, 211; blue, 33 }  ,draw opacity=1 ] (387.99,116.46) .. controls (387.99,116.46) and (387.99,116.46) .. (387.99,116.46) .. controls (387.99,110.82) and (364.88,106.24) .. (336.37,106.24) .. controls (307.86,106.24) and (284.75,110.82) .. (284.75,116.46) ;  
\draw  [draw opacity=0] (511.87,116.46) .. controls (511.87,116.46) and (511.87,116.46) .. (511.87,116.46) .. controls (511.87,116.46) and (511.87,116.46) .. (511.87,116.46) .. controls (511.87,110.82) and (484.14,106.24) .. (449.93,106.24) .. controls (415.72,106.24) and (387.99,110.82) .. (387.99,116.46) -- (449.93,116.46) -- cycle ; \draw  [color={rgb, 255:red, 126; green, 211; blue, 33 }  ,draw opacity=1 ] (511.87,116.46) .. controls (511.87,116.46) and (511.87,116.46) .. (511.87,116.46) .. controls (511.87,116.46) and (511.87,116.46) .. (511.87,116.46) .. controls (511.87,110.82) and (484.14,106.24) .. (449.93,106.24) .. controls (415.72,106.24) and (387.99,110.82) .. (387.99,116.46) ;  
\draw  [draw opacity=0] (532.44,109.32) .. controls (519.71,111.16) and (511.87,113.68) .. (511.87,116.46) -- (584.13,116.46) -- cycle ; \draw  [color={rgb, 255:red, 126; green, 211; blue, 33 }  ,draw opacity=1 ] (532.44,109.32) .. controls (519.71,111.16) and (511.87,113.68) .. (511.87,116.46) ;  

\draw (543.67,158.28) node [anchor=north west][inner sep=0.75pt]    {$\textcolor[rgb]{0.29,0.56,0.89}{G}\textcolor[rgb]{0.29,0.56,0.89}{_{1}}$};
\draw (538.73,64.33) node [anchor=north west][inner sep=0.75pt]    {$\textcolor[rgb]{0.49,0.83,0.13}{G}\textcolor[rgb]{0.49,0.83,0.13}{_{2}}$};

\end{tikzpicture}
\caption{Union $G_1,G_2$}
\label{u2d}
\end{figure}
$G_1$ and $G_2$ are doubling graphs but let us show that $G=G_1 \cup G_2$ is not a doubling graph.

$G_1\cap G_2=\{(\pm x_i,y):x_0=0 \text{  and  } x_i=x_{i-1}+i \text{  with  } i \in \NN\}$. Fix $B$ a ball of $G$ of radius $2k$ centered in $(x_{4k},0)$, see figure \ref{bu2}. Define the set $A=\cup_{i=1}^{k} (x_{5k-i},-k-i)$. For each point in $A$ we can see that $\delta_G\big((x_{5k-i},-k-i),(x_{4k},0)\big)=k+i+k-i=2k$ then $A\subset B$. But for each point of $A$ we have that

  \begin{equation*}
    \begin{split}
    &\delta_G\big((x_{5k-i},-k-i),(x_{5k-j},-k-j)\big)\\&
 \qquad =\min \{\delta_{G_1}\big((x_{5k-i},-k-i),(x_{5k-j},-k-j)\big),\\&
 \qquad \ \delta_{G_1}((x_{5k-i},-k-i), (x_{5k-i},0)+\delta_{G_2}((x_{5k-i},0),(x_{5k-j},0)+\delta_{G_1}\big((x_{5k-j},0),(x_{5k-j},-k-j)\big)\}\\&
 \qquad >k.
    \end{split}
    \end{equation*}
\begin{figure}[h]
\centering
\input{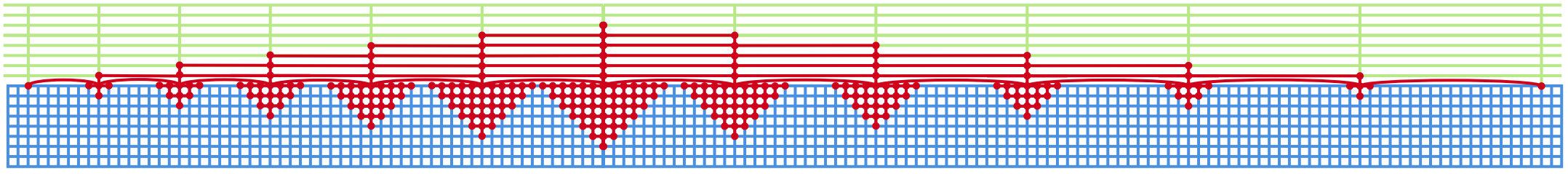}
\caption{Ball with radius 6 of $G$, $B_{\delta_G}(x_{12},6)$.}
\label{bu2}
\end{figure}
    Then, at least $k$ balls of radius $k/2$ are needed to cover $B$. 

Thus, as $k$ increases, the number of balls of radius $k/2$ needed to cover $B$ also increases arbitrarily. Hence, $G$ is not a doubling graph.
\end{example}

\section{Union of doubling graphs with arbitrary intersection}

We have shown that if the intersections are finite sets, then their union is doubling. In this section, we provide conditions under which the union remains doubling even when the intersections between the graphs are infinite sets.

\begin{definition}
Let $G = \bigcup_{i \in I} G_i$. We say that $G$ is compatible with the metrics of the graphs $\{G_i\}_{i\in I}$ if, for every $i \in I$ and for all $v_1, v_2 \in G_i$, we have
$$
\delta_{G_i}(v_1,v_2) = \delta_G(v_1,v_2).
$$
\end{definition}

\begin{theorem}
Let $\{G_i\}_{i=1}^N$ be a finite family of doubling graphs such that $G=\bigcup_{i =1}^N G_i$ is compatible with the metrics of  $\{G_i\}_{i=1}^N$. Then $G$ is a doubling graph.
\label{p1}
\end{theorem}

\begin{proof}
Let $\{G_i\}_{i=1}^N$ be a family of doubling graphs such that $G=\bigcup_{i=1}^N G_i$ is compatible with the metrics of $\{G_i\}_{i=1}^N$. Let $\{C_i\}_{i=1}^N$ be their doubling constants. Let be $I=\{1, \dots, N\}$.

Take a ball of $G$ with radius $2k$.

If $B\cap X_j \neq\emptyset$ and $B\cap \big(\cup_{i\in I, i \neq j} G_i\big)=\emptyset$ then $B$ can be covered by $C_j$ balls of radius $k$.

Let $\hat{I}=\{i_1,\dots, i_k\}$ be a subset of $I$. If $B\cap G_{i_j} \neq\emptyset$ for all $i_j\in \hat{I}$ and $B\cap \big(\cup_{i\in I, i \notin \hat{I}} G_i \big) = \emptyset$. For each $i_j\in \hat{I}$ fix $x_j \in B \cap G_{i_j}$. Since $G$ is compatible with the metrics, then, $B \cap G_{i_j} \subset B_{G_{i_j}}(x_j,4k)$. So, by Lemma \ref{l:doubling}, for each $i_j\in \hat{I}$, $B \cap G_{i_j}$ can be covered by $C_{i_j}^2$ balls of radius $k$. Hence, $B$ can be covered by at most $\sum_{i_j\in \hat{I}}C_{i_j}^2$ balls of radius $k$.

Particularly, if $B\cap G_{i} \neq\emptyset$ for all $i\in I$. Then, $B$ can be covered by at most $\sum_{i\in I}C_{i}^2$ balls of radius $k$.

So for each $k>0$ a ball of $G$ with radius $2k$ can be covered by at most $\sum_{i\in I}C_{i}^2$ balls of radius $k$. Hence, $G$ is a doubling graph.

\end{proof}

For infinite unions of graphs, not even the previous property ensures that the doubling condition is preserved, as illustrated by the Cantor tree example discussed in the introduction, or by the following example.

\begin{example}
Let $\{G_i\}_{i\in \NN\cup \{0\}}$ be a family of doubling graphs given by $G_i=(V_i,E_i)$ where $V_i=\{(i,y): y\in \NN\cup \{0\}\}$ and $E_i=\{\big((i,y),(i,z)\big)\in V_i\times V_i: |y-z|=1\}$ for each $i\in \NN$, and $G_0=(V_0,E_0)$ with $V_0=\{(x,0):x\in \NN\}$ and $E_0=\{\big((x,0),(x',0)\big):|x-x'|=1\}$.

Now, let us define $G=\bigcup_{i \in \NN \cup \{0\}} G_i$, see figure \ref{pei} . For each $i$ we have that for all $v_1,v_2 \in G_i$ $\delta_G(v_1,v_2)=\delta_{G_i}(v_1,v_2)$. But $G$ is not a doubling graph.
\begin{figure}[h]
\centering
\input{peine}
\caption{$G=\bigcup_{i \in \NN \cup \{0\}} G_i$}
\label{pei}
\end{figure}
\end{example}

In fact, the compatibility requirement can be weakened by introducing the concept of quasi-compatibility.

\begin{definition}
Let $G = \bigcup_{i \in I} G_i$. We say that \(G\) is quasi-compatible with the metrics of the graphs $\{G_i\}_{i\in I}$ if, for every $i \in I$ there exists a constant $k_i\geq 1$ such that  for all $v_1, v_2 \in G_i$, 
$$
\delta_{G_i}(v_1,v_2) < k_i \delta_G(v_1,v_2).
$$
\end{definition}

\begin{theorem}
Let $\{G_i\}_{i=1}^N$ be a finite family of doubling graphs such that $G=\bigcup_{i=1}^N G_i$ is quasi-compatible with the metrics of  $\{G_i\}_{i=1}^N$. Then $G$ is a doubling graph.
\label{p2}
\end{theorem}

\begin{proof}
First, let us consider a family doubling graphs $\{G_i\}_{i=1}^N$ such that $G=\bigcup_{i=1}^N G_i$  is quasi-compatible with their metrics. Let $\{C_i\}_{i=1}^N$ be their doubling constants and $I=\{1, \dots, N\}$.

Since $G$ is quasi-compatible with the metrics of  $\{G_i\}_{i=1}^N$, we can consider the family of constants $\{k_i\}_{i=1}^N$ such that, for each $i\in I$,
$$
\delta_{G_i}(x,y)<{k_i}\delta_G(x,y)
$$
for all $x,y\in G_i$.

For each one of these constants we can choose $k_i'$ such that $k_i \le 2^{k_i'}$ to simplify the calculations.

Let us take $B$ a ball of $G$ with radius $2R$.

If $B\cap X_j \neq\emptyset$ and $B\cap \big(\bigcup_{i=1,i\neq j}^N G_i\big)=\emptyset$, then $B$ can be covered by $C_j$ balls of radius $R$.

Otherwise, let $\hat{I}=\{i_1,\dots, i_k\}$ be the subset of $I$ such that for all $i_j\in \hat{I}$, $B\cap G_{i_j}\neq \emptyset$ and $B\cap \big(\bigcup_{i\in I, i \notin \hat{I}} G_i\big)=\emptyset$. For each $i_j \in \hat{I}$ let us fix $x_j\in B\cap G_{i_j}$. Then, for all $y\in G_{i_j}\cap B$, we have,
$$
\delta_{G_{i_j}}(x_j,y)<2^{k_{i_j}'}\delta_G(x_j,y)<2^{k_{i_j}'}4R.
$$
So, 
$$
G_{i_j}\cap B\subset B_{i_j}(x_j,2^{k_{i_j}'}4R)=\{y\in G_{i_j}: \delta_{G_{i_j}}(x_j,y)<2^{k_{i_j}'}4R\}.
$$
Hence, by Lemma \ref{l:doubling}, $G_{i_j}\cap B$ can be covered by $C_{i_j}^{2+k_{i_j}'}$ balls of radius $R$.

It follows that $B$ can be covered by $\sum_{i_j\in \hat{I}} C_{i_j}^{2+k_{i_j}'}$ balls of radius $R$.

Hence, any ball of $G$ with radius $2R$ can be covered with at most $\sum_{i=1}^N C_i^{2+k_i'}$ balls of radius $R$, so $G$ is a doubling graph. Hence $G$ is a doubling graph.
\end{proof}

Although this condition is quite general, it does not cover all cases in which gluing two doubling graphs yields a doubling graph as the following example shows.

\begin{example}\label{exp} Consider the graphs $G_1$ and $G_2$ defined as,
$$\begin{aligned} 
V_1=&\{(x,y)\in \ZZ^2: y=0 \text{ and }x \geq 0 \text{ or } x=0 \text{ and } y\geq 0\}\cup \{(x,x)\in \ZZ^2: x < 0\},\\
E_1=&\{\big((x,0),(-x,-x)\big)\in V_1 \times V_1 : x\geq 0\}\cup \{\big((0,y),(-y,-y)\big)\in V_1 \times V_1 : y\geq 0\}\\
&\cup \{\big((x,x),(x-1,x-1)\big)\in V_1 \times V_1 : x\leq 0\}, \\
V_2=&\{(x,y)\in \ZZ^2: y\geq 0 \text{ and }x \geq 0\},\\
E_2=&\{\big((x,y),(z,t)\big)\in V_2 \times V_2 : y=t \text{ and }x =z+1 \text{ or } x=z \text{ and } y=t+1\}.
\end{aligned}$$
See Figure \ref{B}. These graphs are doubling and gluing them along their common vertices yields a graph that is also doubling. However, notice that for every integer $x>0$, $\delta_{G_2}\big((x,0),(0,x)\big)=2x$ while $\delta_{G}\big((x,0),(0,x)\big)=2$. Therefore, there is no constant $k_2\geq 0$ such that the inequality 
$$\delta_{G_2}(v_1,v_2)\leq k_2 \delta_G(v_1,v_2)$$
holds for every pair of vertices $v_1,v_2\in V_2$. Then, the metrics are not quasi-compatible.
\begin{figure}[h]
\centering
\input{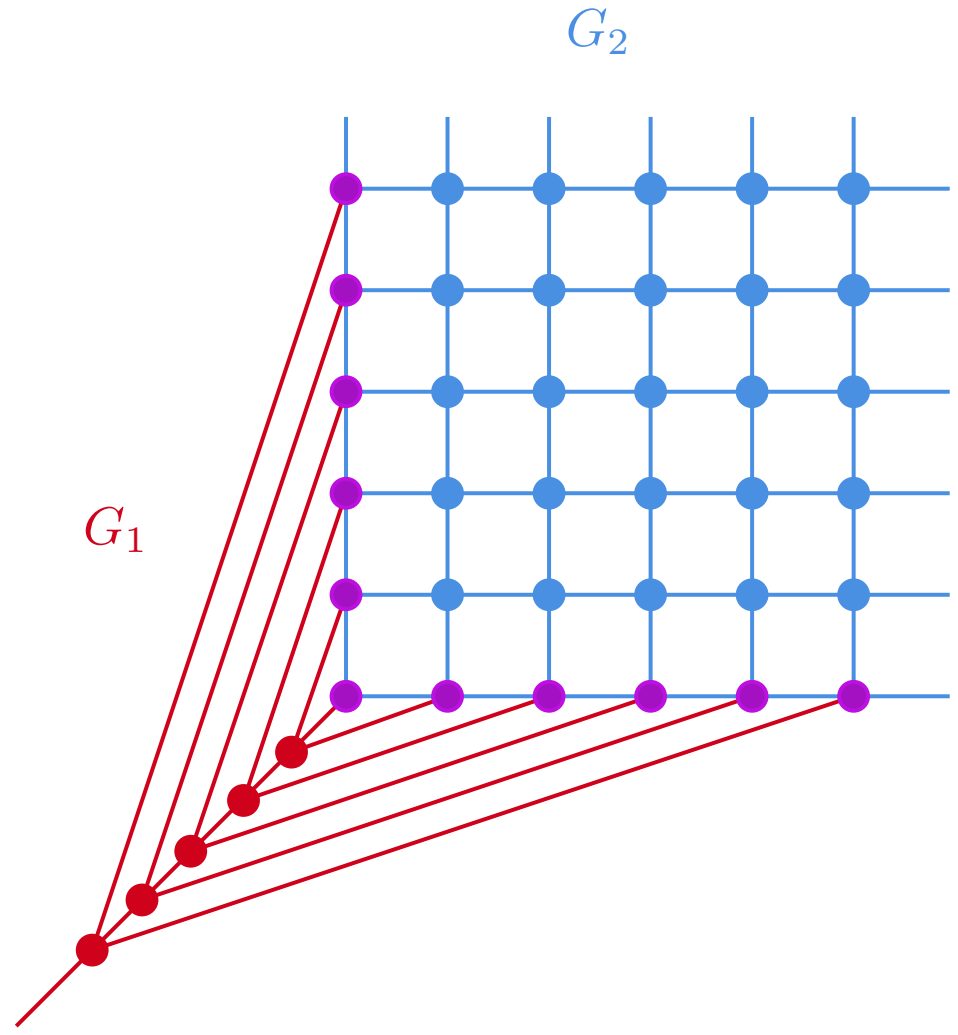}
\caption{$G= G_1 \cup G_2$}
\label{B}
\end{figure}
\end{example}

\section{Decomposition of doubling graphs}\label{s:decomposition}

In this section, we study the case in which the graph obtained by gluing a family of graphs is doubling. Does this imply that each graph in the family is doubling?

In general, even if the resulting graph is doubling, this does not imply that the glued subgraphs are doubling. Even when gluing a finite number of non-doubling graphs, the resulting graph can still be doubling.

\begin{example}
Consider the graphs $G_1=(V_1,E_1)$ and $G_2=(V_2,E_2)$ such that $V_1 = \{(x,y) \in \mathbb{Z}^2 : x,y \ge 0\}$, $E_1 = \{\big((x,y),(z,y)\big) \in V_1 \times V_1 : |x-z|=1\} \cup \{\big((0,y),(0,z)\big) \in V_1 \times V_1 : |y-z|=1\}$, and $V_2 = \{(x,y) \in \mathbb{Z}^2 : x,y \ge 0\}$, $E_2 = \{\big((x,y),(x,z)\big) \in V_2 \times V_2 : |y-z|=1\} \cup \{\big((x,0),(z,0)\big) \in V_2 \times V_2 : |x-z|=1\}$. Neither graph is doubling, but we can glue them to obtain the graph $G=(V,E)$, where $V = \{(x,y) \in \mathbb{Z}^2 : x,y \ge 0\}$ and $E = \{\big((x,y),(t,z)\big) \in V \times V : (x=t \text{ and } |y-z|=1) \text{ or } (y=z \text{ and } |x-t|=1)\}$, which is indeed a doubling graph.
\end{example}

\begin{theorem}\label{t:decomp}
Let $\{G_i\}_{i\in I}$ be a (possibly infinite) family of graphs such that $G=\cup_{i\in I} G_i$ is quasi-compatible with the metrics of $\{G_i\}_{i\in I}$. If $G$ is a doubling graph, then each $G_i$ is a doubling graph. 
\end{theorem}

\begin{proof}
Suppose $G$ is doubling with constant $C$ and fix any $G_i$. Then, by definition, there exists a constant $k_i\geq 1$ such that for every $x,y \in G_i$,
$$\delta_{G_i}(x,y)\leq k_i\delta_G(x,y).$$
Let us consider $l_i\in\ZZ$ a constant such that $k_i\leq 2^{l_i}$.

Consider any $x\in G_i$ and any ball $B_{G_i}(x,2k)$. Then $B_{G_i}(x,2k)\subset B_G(x,2^{l_i+1}k)$ and, by Lemma \ref{l:doubling}, this ball can be covered by $C^{2l_i+3}$ balls of radius $\frac{k}{2^{l_i+1}}$. Let $\{B_j\}_{j=1}^{C^{2l_i+3}}$ be these balls and define the set $J=\{j\in \{1,\dots, C^{2l_i+3}\}: G_i\cap B_j \neq \emptyset\}$.

For each $j\in J$ let us fix $x_j\in G_i\cap B_j$. Then for all $y \in G_i\cap B_j$ we have that,
$$\delta_{G_i}(x_j,y)\leq k_i \delta_G(x_j,y)\leq \frac{2k k_i}{2^{l_i+1}}\leq \frac{k k_i}{k_i}=k.$$

From this we get,

$$B_{G_i}(x,2k)\subset \cup_{j\in J} B_j \cap G_i \subset \cup_{j\in J} B_{G_i}(x_j,k).$$

Thus, in $(G_i,\delta_{G_i})$ every ball of radius $2k$  can be covered by at most $C^{2l_i+3}$ balls of radius $k$ and $G_i$ is a doubling graph.
\end{proof}

Example \ref{exp} also shows that the condition above is not necessary.

\section{Gluing and decomposing doubling metric spaces}

To work more generally with metric spaces, we must be careful when performing gluing operations. If we attempt to glue simply by identifying subsets of both metric spaces as we did in the case of graphs, we may find that the resulting space fails to be a metric space. See \cite{BH}. In this article, the gluing of metric spaces that we shall use is the $L$-bi-Lipschitz gluing.

\begin{definition}[$L$-bi-Lipschitz Gluing Metric]
Given two metric spaces, $(X_1, d_{X_1})$ and $(X_2, d_{X_2})$, consider two closed subspaces, $U \subset X_1$ and $V \subset X_2$, and an $L$-bi-Lipschitz bijection between them; that is, a bijection $f \colon U \to V$ satisfying
$$
L^{-1} d_{X_1}(a, b) \le d_{X_2}(f(a), f(b)) \le L \, d_{X_1}(a, b) \quad \text{for all } a, b \in U,
$$
for some constant $L\geq 1$.The space resulting from gluing along this map is
$$
X := (X_1 \sqcup X_2) / \sim
$$
where $\sim$ is the equivalence relation generated by the identifications $a \sim f(a)$ for all $a \in U$. 
\end{definition}
\begin{proposition}\cite[Proposition 3.2.8]{BR24}
$X$ is a metric space equipped with the metric $d_X$ defined for any two equivalence classes $[a], [b] \in X$ by:
\begin{enumerate}
    \item For points in different spaces ($a \in X_1$ and $b \in X_2$):
    $$
    d_X([a], [b]) = \inf_{u \in U} \left\{ d_{X_1}(a, u) + d_{X_2}(f(u), b) \right\}.
    $$

    \item For points in $X_1$ ($a, b \in X_1$):
    $$
    d_X([a], [b]) = \min \left\{ d_{X_1}(a, b), \, \inf_{u_1, u_2 \in U} \left\{ d_{X_1}(a, u_1) + d_{X_2}\big(f(u_1), f(u_2)\big) + d_{X_1}(u_2, b) \right\} \right\}.
    $$

    \item For points in $X_2$ ($a, b \in X_2$):
    $$
    d_X([a], [b]) = \min \left\{ d_{X_2}(a, b), \, \inf_{v_1, v_2 \in V} \left\{ d_{X_2}(a, v_1) + d_{X_1}\big(f^{-1}(v_1), f^{-1}(v_2)\big) + d_{X_2}(v_2, b) \right\} \right\}.
    $$
\end{enumerate}
\end{proposition}


\begin{obs}
Let $X=X_1 \cup_f X_2$ be the metric space given by gluing $X_1$ and $X_2$ along the $L$-bi-Lipschitz bijection $f:U \to V$, where $U\subset X_1$ and $V \subset X_2$ are closed subspaces. Then, for all $x_1,x_2\in X_i$ we have:

$$d_X([x_1],[x_2]) \leq d_{X_i}(x_1,x_2)\leq  L d_X([x_1],[x_2]).$$
\end{obs}

\begin{theorem}\label{th:metric}
Let $X=X_1 \cup_f X_2$ be the metric space given by gluing two metric spaces, $X_1$ and $X_2$, along the $L$-bi-Lipschitz bijection $f:U \to V$, where $U\subset X_1$ and $V \subset X_2$ are closed subspaces. Then $X$ is a doubling metric space if and only if $X_1$ and $X_2$ are metric doubling spaces.
\end{theorem}

\begin{proof}
Let us consider $N\in\ZZ$ a constant such that $L\leq 2^{N}$.

We first prove the forward implication. Let us suppose that $X=X_1\cup_f X_2$ is a doubling metric space and $C$ its doubling constant. Fix $i \in \{1,2\}$ and a point $x \in X_i$. Given $k\geq 0$, we want to show that the ball $B_{d_{X_i}}(x,2k)$ can be covered by at most $C_i$ balls of $X_i$ with radius $k$, where $C_i$ is a constant that does not depend on $x$. 

Since $X$ is doubling, by Lemma \ref{l:doubling}, $B_{d_X}([x],2k)$ can be covered by $C^{N+2}$ balls, $B_1,\dots,B_{C^{N+2}}$, of radius $\frac{2k}{2^{N+2}}$. Define $J=\{j\in \{1,\dots, C^{N+2}\}: \textrm{ exists } y \in X_i \textrm{ such that } [y] \in B_j\}$. Now, for each $j\in J$ let us take $y_j\in X_i$ such that $[y_j]\in B_j$. Then, for all $z\in X_i$ such that $[z]\in B_j$, we have that,

$$d_{X_i}(z,y_j)\leq L d_X([z],[y_j])\leq \frac{2^{N+1}k}{2^{N+1}}=k.$$

So, $z\in B_{d_{X_i}}(y_j,k)$ and it follows that,
 $$B_{d_{X_i}}(x,2k)\subset \cup_{j\in J} B_{d_{X_i}}(y_j,k).$$
 Then, $B_{d_{X_i}}(x,2k)$ can be covered by at most $C^{N+2}$ balls of radius $k$ in $X_i$. Hence, $X_i$ is a doubling metric space.
 
We now prove the reverse implication. Fix $[x]\in X$ and $k\geq 0$. Consider the ball $B=B_{d_X}([x],2k)$ and let us defined the set $S=\{y\in X_1\sqcup X_2: [y]\in B\}$.
\begin{itemize}
\item If $S \cap X_1 \neq \emptyset$ and $S\cap X_2=\emptyset$, then, for $x_0\in X_1$ with $[x_0]=[x]$ and for all $y\in S$ we have that,

$$d_{X_1}(x_0,y)\leq 2^{N} d_X([x_0], [y])\leq 2^{N+1}k.$$

Thus, $S\subset B_{d_{X_1}}(x_0,2^{N+1}k)$. Hence, since $X_1$ is doubling and considering $C_1$ its doubling constant, $S$ can be covered by $C_1^{N+1}$ balls of radius $k$ in $X_1$. Let $\{x_j\}_{j=1}^{C_1^{N+1}}$ be the centers of theses balls, then,

$$B\subset \cup_{j=1}^{C^{N+1}} B_{d_X}([x_j],k),$$
so $B$ can be covered by at most $C_1^{N+1}$ balls in $X$.

\item If $S \cap X_2 \neq \emptyset$ and $S\cap X_1=\emptyset$, by the same argument as before, $B$ can be covered by at most $C_2^{N_2+1}$ balls in $X$, with $C_2$ the doubling constant of $X_2$.

\item If $S \cap X_1 \neq \emptyset$ and $S\cap X_2\neq \emptyset$. Then, there must be some point $a\in U$ such that $a\in S$ and $d_X([a],[x])\leq 4k$ for every $[x]\in B$. Therefore,  $S \cap X_1\subset B_{d_{X_1}}(a,2^{N+2}k)$ and $S \cap X_2\subset B_{d_{X_2}}(f(a),2^{N+2}k)$.

Repeating the argument above, $X_1 \cap S$ can be covered by $C_1^{N+2}$ balls of radius $k$ in $X_1$, and $X_2 \cap S$ can be covered by $C_2^{N+2}$ balls of radius $k$ in $X_2$. Let $\{y_j\}_{j=1}^{C_1^{N+2}}$ and $\{z_l\}_{l=1}^{C_2^{N+2}}$ denote the respective sets of centers of these balls. Then, the following holds,

$$B\subset \Big(\cup_{j=1}^{C_1^{N+2}}B_{d_X}([y_j],k)\Big) \cup \Big(\cup_{l=1}^{C_2^{N+2}}B_{d_X}([z_l],k)\Big).$$
 Hence, $B$ can be covered by at most $C_1^{N+2}+C_2^{N+2}$ balls of radius $k$ in $X$.
\end{itemize}

So, we have proved that every ball of radius $2k$ in $X$ can be covered by at most $C_1^{N+2}+C_2^{N+2}$ balls of radius $k$. Therefore, $X$ is a doubling metric space.
\end{proof}

\medskip

\textbf{AI declaration}. This manuscript is part of a PhD thesis. No AI has been used at any stage of this work.

\end{document}